\documentclass{appolb}
\usepackage{bm}
\usepackage{graphicx}
\usepackage{amsfonts}

\title{
Waves along fractal coastlines: From fractal arithmetic to wave equations}
\author{
Marek Czachor\\
Katedra Fizyki Teoretycznej i Informatyki Kwantowej\\ 
Politechnika Gda\'nska \\
ul. G. Narutowicza 11/12, 80-233 Gda\'nsk, Poland}
\begin{document}
\maketitle

\newcommand{\be}{\begin{eqnarray}}
\newcommand{\ee}{\end{eqnarray}}

\begin{abstract}
Beginning with addition and multiplication intrinsic to a Koch-type curve we
formulate and solve wave equation describing wave propagation along a fractal coastline. As opposed to examples known from the literature, we do not replace the fractal by the continuum in which it is embedded. This seems to be the first example of a truly intrinsic description of wave propagation along a fractal curve. The theory is relativistically covariant under an appropriately defined Lorentz group.
\end{abstract}


\section{Non-Newtonian calculus}

Consider two sets $\mathbb{X}$ and $\mathbb{Y}$  whose cardinality is continuum, and a function $A:\mathbb{X}\to\mathbb{Y}$. 
There exist bijections $f_\mathbb{X}$, $f_\mathbb{Y}$, $g_\mathbb{X}$, $g_\mathbb{Y}$, such that the diagram
\be
\begin{array}{rcl}
\mathbb{R}                & \stackrel{\tilde B}{\longrightarrow}   & \mathbb{R}\\
g_\mathbb{X}{\Big\uparrow}   &                                     & {\Big\uparrow}g_\mathbb{Y}   \\
\mathbb{X}                & \stackrel{A}{\longrightarrow}       & \mathbb{Y}               \\
f_\mathbb{X}{\Big\downarrow}   &                                     & {\Big\downarrow}f_\mathbb{Y}   \\
\mathbb{R}                & \stackrel{\tilde A}{\longrightarrow}   & \mathbb{R}
\end{array}\nonumber
\ee
is commutative. The functions $\tilde A$ and $\tilde B$ are defined by the diagram. 
It is natural to think of $\mathbb{X}$ and $\mathbb{Y}$ in terms of one-dimensional manifolds whose global charts are defined by the bijections.

In differential topology and geometry, a derivative of $A:\mathbb{X}\to\mathbb{Y}$ would be a function $A':\mathbb{X}\to\mathbb{Y}$ defined by  $\widetilde{A'}(r)=d\tilde A(r)/dr$. Of course, since
\be
\tilde A = f_\mathbb{Y}\circ g_\mathbb{Y}^{-1}\circ \tilde B\circ g_\mathbb{X}\circ f_\mathbb{X}^{-1}
= \varphi_\mathbb{Y}^{-1}\circ \tilde B\circ \varphi_\mathbb{X},\nonumber
\ee
a derivative of $A$ can be equivalently defined in terms of $\tilde B$, provided $\varphi_\mathbb{X}$ and $\varphi_\mathbb{Y}^{-1}$ are at least $C^1$. A transition between the two forms is determined by the chain rule for derivatives.

In the arithmetic approach to differentiation \cite{GK,G79,G83,MC,ACK1,ACK2,ACK3,MC2,Blaszak} one starts from a different perspective. In the first step, one employs the bijections to turn $\mathbb{X}$ and $\mathbb{Y}$ into fields isomorphic to 
$\mathbb{R}$. Explicitly, one defines the arithmetic operations in $\mathbb{X}$ (addition, subtraction, multiplication, division) by
\be
x\oplus_\mathbb{X} y &=& f_\mathbb{X}^{-1}\big(f_\mathbb{X}(x)+f_\mathbb{X}(y)\big),\\
x\ominus_\mathbb{X} y &=& f_\mathbb{X}^{-1}\big(f_\mathbb{X}(x)-f_\mathbb{X}(y)\big),\\
x\odot_\mathbb{X} y &=& f_\mathbb{X}^{-1}\big(f_\mathbb{X}(x)f_\mathbb{X}(y)\big),\\
x\oslash_\mathbb{X} y &=& f_\mathbb{X}^{-1}\big(f_\mathbb{X}(x)/f_\mathbb{X}(y)\big),
\ee
and analogously in $\mathbb{Y}$. This type of arithmetic is a special case of a general non-Diophantine arithmetic discussed by Burgin \cite{Burgin0,Burgin1,Burgin2,Burgin3}.  The case of a linear $f$ was extensively studied in \cite{Benioff,Benioff1,Benioff2} with emphasis on distinguishing between numbers, treated abstractly, and their representations and values.

The topologies of  $\mathbb{X}$ and $\mathbb{Y}$  are induced by the bijections from the topology of $\mathbb{R}$. 
Let the limit $x\to x_0\in\mathbb{X}$ be defined by the formula
\be
\lim_{x\to x_0}A(x)=f^{-1}_\mathbb{Y}\left(\lim_{r\to f_\mathbb{X}(x_0)}\tilde A(r)\right).
\ee
The derivative of $A$ can be expressed in terms of limits in three equivalent ways
\be
\frac{DA(x)}{Dx}
&=& \lim_{h\to 0}\Big(A\big(x\oplus_{\mathbb{X}}f^{-1}_{\mathbb{X}}(h)\big)\ominus_{\mathbb{Y}}A(x)\Big)\oslash_{\mathbb{Y}}f^{-1}_{\mathbb{Y}}(h)
\\
&=& \lim_{h\to 0}\Big(A\big(x\oplus_{\mathbb{X}}h_{\mathbb{X}}\big)\ominus_{\mathbb{Y}}A(x)\Big)\oslash_{\mathbb{Y}}h_{\mathbb{Y}}
\\
&=& \lim_{h\to 0_{\mathbb{X}}}\Big(A\big(x\oplus_{\mathbb{X}}h\big)\ominus_{\mathbb{Y}}A(x)\Big)\oslash_{\mathbb{Y}}f(h),
\ee
where $f=f^{-1}_\mathbb{Y}\circ f_\mathbb{X}$. Here $0_{\mathbb{X}}$ is the neutral element of addition in $\mathbb{X}$.
This type of derivative was investigated in a systematic way for the first time in \cite{GK}, for the case where $\mathbb{X}$ and 
$\mathbb{Y}$ were subsets of $\mathbb{R}$, while $f_\mathbb{Y}$ and $ f_\mathbb{X}$ were continuous in the metric topology of $\mathbb{R}$. The derivative was rediscovered by myself in a fractal context \cite{MC}. The main difference between the formalism from \cite{GK} and my approach is that now the derivative is applicable to all sets whose cardinality is continuum, such as Sierpi{\'n}ski-type fractals \cite{ACK3}, which obviously do not have to be subsets of $\mathbb{R}$, and in typical examples the bijections are discontinuous in metric topologies of $\mathbb{X}$ and $\mathbb{Y}$. This counterintuitive possibility opened by the non-Newtonian calculus is especially useful in fractal applications. Just to give one example, a construction of Fourier transforms on arbitrary Cantor sets is in the non-Newtonian framework basically trivial \cite{ACK2}, simultaneously circumventing various impossibility theorems known from the more traditional approach \cite{JP,PJ}. The arithmetic perspective is simultaneously applicable to all the other aspects of mathematical modeling, including algebraic or probabilistic methods. The freedom of choice of arithmetic plays a role of a universal symmetry of any mathematical model.

An equivalent and very convenient form of the derivative is
\be
\frac{DA(x)}{Dx}
&=&
f^{-1}_{\mathbb{Y}}\left(
\frac{d\tilde A\big(f_{\mathbb{X}}(x)\big)}{df_{\mathbb{X}}(x)}
\right).
\label{def non-N D}
\ee
The derivative is Newtonian if $\mathbb{X}$ and $\mathbb{Y}$ are subsets of $\mathbb{R}$, and $f_{\mathbb{X}}(x)=x$,$f_{\mathbb{Y}}(y)=y$ are the identity maps. If the bijections are less trivial, one speaks of non-Newtonian derivatives. 

Of particular interest is the non-Newtonian version of the chain rule. Consider the diagram
\be
\begin{array}{rrrrr}
\mathbb{X}                & \stackrel{A}{\longrightarrow}       & \mathbb{Y}              & \stackrel{B}{\longrightarrow}       & \mathbb{Z}   \\
f_\mathbb{X}{\Big\downarrow}   &                                     & f_\mathbb{Y}{\Big\downarrow}   &                                     & f_\mathbb{Z}{\Big\downarrow}   \\
\mathbb{R}                & \stackrel{\tilde A}{\longrightarrow}   & \mathbb{R} & \stackrel{\tilde B}{\longrightarrow}   & \mathbb{R}
\end{array}\nonumber
\ee
Then
\be
\frac{D(B\circ A)(x)}{Dx}
=
f^{-1}_\mathbb{Z}\left[f_\mathbb{Z}\Bigg(\frac{D B\big(A(x)\big)}{DA(x)}\Bigg) f_\mathbb{Y}\Bigg(\frac{D A(x)}{Dx}\Bigg)\right].
\label{chain1}
\ee
For a composition of three functions,
\be
\mathbb{W} \stackrel{A}{\longrightarrow} \mathbb{X} \stackrel{B}{\longrightarrow}  \mathbb{Y}   \stackrel{C}{\longrightarrow}  \mathbb{Z},
\label{chain3}
\ee
one finds
\be
\frac{D C\circ B\circ A(x)}{Dx}
=
f^{-1}_\mathbb{Z}\left[f_\mathbb{Z}\Bigg(\frac{D C\big[B\big(A(x)\big)\big]}{DB\big(A(x)\big)}\Bigg) 
f_\mathbb{Y}\Bigg(\frac{D B\big(A(x)\big)}{DA(x)}\Bigg) f_\mathbb{X}\Bigg(\frac{D A(x)}{Dx}\Bigg)\right].
\label{cor 3}
\ee
The latter case is important since it allows us to better understand the structure of the non-Newtonian derivative. Indeed, let the three functions be the ones occurring in the definition of $\mathbb{X} \stackrel{A}{\longrightarrow} \mathbb{Y}$, i.e.
\be
\mathbb{X} \stackrel{f_\mathbb{X}}{\longrightarrow} \mathbb{R} \stackrel{\tilde A}{\longrightarrow}  \mathbb{R}   \stackrel{f_\mathbb{Y}^{-1}}{\longrightarrow}  \mathbb{Y}.
\ee
Now, directly from definition one checks that
\be
\frac{Df_\mathbb{X}(x)}{Dx}=1,\quad
\frac{Df_\mathbb{Y}(x)}{Dx}
=1,\quad
\frac{Df_\mathbb{X}^{-1}(x)}{Dx}
=1_\mathbb{X},\quad
\frac{Df_\mathbb{Y}^{-1}(x)}{Dx}
=
1_\mathbb{Y}.\label{4 1}
\ee
The chain rule implies
\be
\frac{DA(x)}{Dx}
&=&
\frac{D (f_\mathbb{Y}^{-1}\circ \tilde A\circ f_\mathbb{X})(x)}{Dx}\nonumber\\
&=&
f^{-1}_\mathbb{Y}\left[f_\mathbb{Y}\Bigg(\frac{D f_\mathbb{Y}^{-1}\big[\tilde A\big(f_\mathbb{X}(x)\big]}{D\tilde A\big(f_\mathbb{X}(x)\big)}\Bigg) 
f_\mathbb{R}\Bigg(\frac{D \tilde A\big(f_\mathbb{X}(x)\big)}{Df_\mathbb{X}(x)}\Bigg) f_\mathbb{R}\Bigg(\frac{D f_\mathbb{X}(x)}{Dx}\Bigg)\right].
\nonumber
\ee
The arithmetic in $\mathbb{R}$ is Diophantine, $f_\mathbb{R}(x)=x$, and thus 
\be
\frac{D \tilde A\big(f_\mathbb{X}(x)\big)}{Df_\mathbb{X}(x)}
&=&
\frac{d \tilde A\big(f_\mathbb{X}(x)\big)}{df_\mathbb{X}(x)}\nonumber
\ee
is Newtonian. Derivatives (\ref{4 1}) imply
\be
\frac{DA(x)}{Dx}
=
f^{-1}_\mathbb{Y}\left(\frac{d \tilde A\big(f_\mathbb{X}(x)\big)}{df_\mathbb{X}(x)}\right),\nonumber
\ee
and we reconstruct our definition of the derivative. One concludes that the bijections behave as identity maps with respect to non-Newtonian derivatives they define, no matter how weird the bijections themselves actually are.

The integral is defined in a way guaranteeing the fundamental laws of calculus, relating derivatives and integrals:
\be
\int_Y^X A(x)Dx
&=&
f^{-1}_{\mathbb{Y}}\left(\int_{f_{\mathbb{X}}(Y)}^{f_{\mathbb{X}}(X)}\tilde A(x)dx\right)
\label{int}
\ee
where $\int \tilde A(x)dx$ is the usual (say, Lebesgue) integral of a function $\tilde A:\mathbb{R}\to \mathbb{R}$.
One proves that
\be
\frac{D}{DX}\int_Y^X A(x)Dx
&=&A(X),\label{int1}\\
\int_Y^X\frac{DA(x)}{Dx}Dx
&=&
A(X)\ominus_{\mathbb{Y}} A(Y)\label{int2}.
\ee

Let us now see how it works in the simple but instructive case of $f(x)=x^3$. The manifold in question is $\mathbb{X}=\mathbb{R}$. Let the two (global) charts be given by $f(x)=x^3$ and $g(x)=x$. Their composition $g\circ f^{-1}$ is not a diffeomorphism if the differentiation is understood in the Newtonian way. Apparently, $f(x)=x^3$ does not define a differentiable structure on $\mathbb{R}$. In the standard Newtonian formalism the only structure we have at our disposal is $C^0$.

The arithmetic approach begins with arithmetic operations intrinsic to $\mathbb{X}$,
\be
x\oplus y &=& f^{-1}\big(f(x)+f(y)\big)=\sqrt[3]{x^3+y^3},\label{+}\\
x\ominus y &=& f^{-1}\big(f(x)-f(y)\big)=\sqrt[3]{x^3-y^3},\label{-}\\
x\odot y &=& f^{-1}\big(f(x)f(y)\big)=\sqrt[3]{x^3y^3}=xy,\label{.}\\
x\oslash y &=& f^{-1}\big(f(x)/f(y)\big)=\sqrt[3]{x^3/y^3}=x/y.\label{/}
\ee
Let us stress again that $f$ is, by construction, a field isomorphism of $(\mathbb{R},+,\cdot)$ and $(\mathbb{R},\oplus,\odot)$. Therefore, $\oplus$ and $\odot$ are commutative and associative, and $\odot$ is distributive with respect to $\oplus$. 
The neutral elements of $\oplus$ and $\odot$, $0'$ and $1'$, are the standard ones: $0'=f^{-1}(0)=\sqrt[3]{0}=0$, $1'=f^{-1}(1)=\sqrt[3]{1}=1$. 
Although multiplication is unchanged, the link between addition and multiplication is a subtle one, as can be seen in the following example
\be
x\oplus\dots \oplus x
&=&
\sqrt[3]{x^3+\dots+x^3}\quad (\textrm{$n$ times})\\
&=&
\sqrt[3]{n}x =f^{-1}(n)x.
\ee
The inverse bijection $f^{-1}(x)=\sqrt[3]{x}$ is continuous but not Newtonian differentiable at $x=0$, hence the loss of the Newtonian diffeomorphism property. Still, the derivative of a function $A:\mathbb{X}\to \mathbb{X}$,
\be
\frac{D A(x)}{Dx}
&=&
\lim_{h\to 0}\Big(A(x\oplus h)\ominus A(x)\Big)\oslash h\label{der0},
\ee
is well defined. The non-Newtonian $D/Dx$ satisfies all the basic rules of differentiation, of course with respect to the new arithmetic:

\noindent
(a) The Leibniz rule.
\be
\frac{D A(x)B(x)}{Dx}
&=&
\lim_{h\to 0}\Big(A(x\oplus h)B(x\oplus h)\ominus A(x)B(x)\Big)\oslash h\\
&=&
\sqrt[3]{\left(\frac{DA(x)}{Dx}B(x)\right)^3 +\left(A(x)\frac{DB(x)}{Dx}\right)^3}\nonumber\\
&=&
\frac{D A(x)}{Dx}B(x)\oplus A(x)\frac{D B(x)}{Dx}.
\ee
(b) Linearity.
\be
\frac{D A(x)\oplus B(x)}{Dx}
&=&
\lim_{h\to 0}\Big(A(x\oplus h)\oplus B(x\oplus h)\ominus \big(A(x)\oplus B(x)\big)\Big)\oslash h\\
&=&
\sqrt[3]{\left(\frac{DA(x)}{Dx}\right)^3+ \left(\frac{DB(x)}{Dx}\right)^3}\\
&=&
\frac{DA(x)}{Dx} \oplus\frac{DB(x)}{Dx}.
\ee
(c) The chain rule.
Denoting
\be
H &=& B(x\oplus h)\ominus B(x),
\ee
we obtain
\be
\frac{D A\big(B(x)\big)}{Dx}
&=&
\lim_{h\to 0}\sqrt[3]{\frac{\Big(A\big(B(x)\oplus H\big)\Big)^3- \Big(A\big(B(x)\big)\Big)^3}{h^3}}\\
&=&
\lim_{H\to 0}\frac{A\big(B(x)\oplus H\big)\ominus A\big(B(x)\big)}{H}
\lim_{h\to 0}\frac{B(x\oplus h)\ominus B(x)}{h}\\
&=&
\frac{DA\big(B(x)\big)}{DB(x)}\frac{DB(x)}{Dx}.
\ee

\begin{figure}
\includegraphics{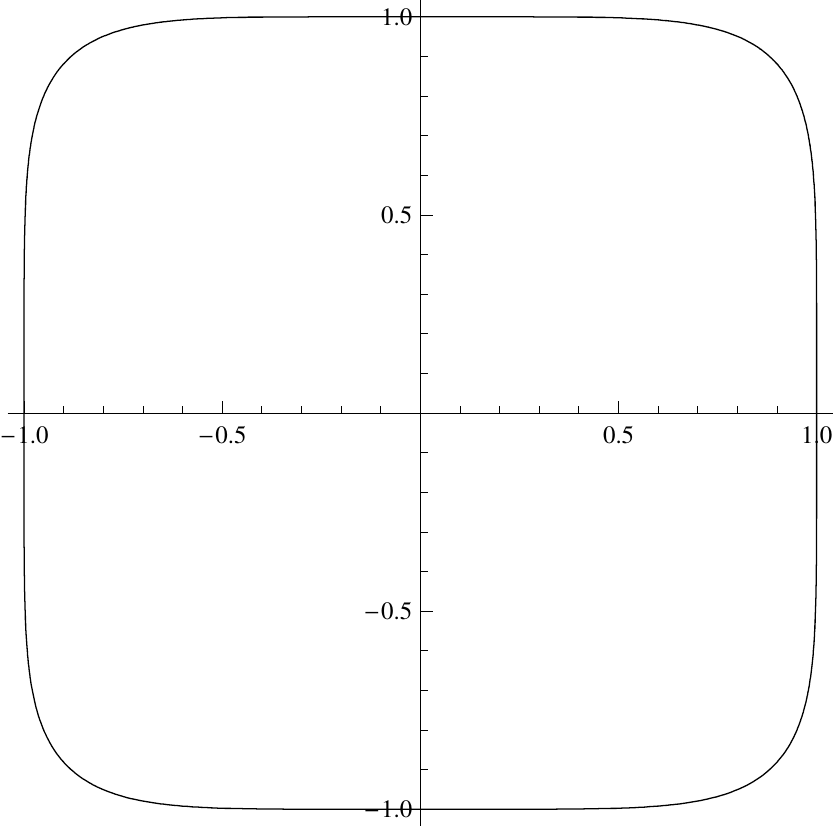}
\caption{The circle $x\mapsto (\textrm{Cos }x,\textrm{Sin }x)$, $0\leq x<(2\pi)^{1/3}$, with trigonometric functions given by 
(\ref{sin x3})--(\ref{cos x3}).}
\label{Fig1}
\end{figure}

\begin{figure}
\includegraphics{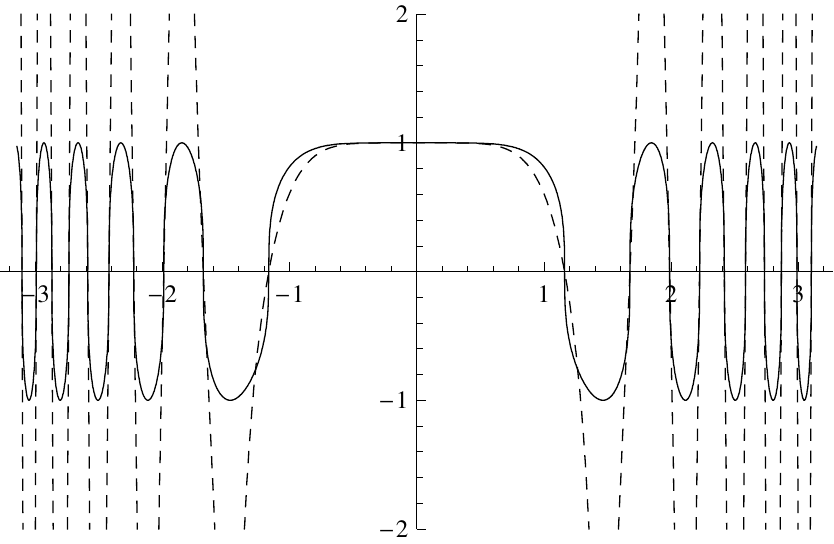}
\caption{The non-Newtonian derivative $D\textrm{Sin }x/Dx=\textrm{Cos }x$ (full, Eq. (\ref{DSin})),  as compared to the standard Newtonian $d\textrm{Sin }x/dx$ (dashed, Eq.  (\ref{dSin})). The singular behavior of the dashed curve follows from Newtonian non-differentiability of $f^{-1}(x)=\sqrt[3]{x}$ at $x=0$. In contrast, the non-Newtonian derivative is non-singular since  $f$ and $f^{-1}$ get differentiated in a non-Newtonian way, yielding trivial derivatives.}
\label{Fig2}
\end{figure}

The non-Newtonian derivate has interesting implications for differential equations. For example, the unique solution of
\be
\frac{D A(x)}{Dx}=A(x),\quad A(0)=1,\label{exp eq}
\ee
is 
\be
A(x)=e^{x^3/3}=f^{-1}\left(e^{f(x)}\right),\label{exp sol}
\ee
as one can verify directly from definition (\ref{der0}).
The exponent satisfies the usual law
\be
A(x_1\oplus x_2)
&=&
e^{(x_1^3+x_2^3)/3}
=A(x_1)\odot A(x_2).
\ee
One can similarly verify that 
\be
\textrm{Sin }x &=& \sqrt[3]{\sin (x^3)},\label{sin x3}\\
\textrm{Cos }x &=& \sqrt[3]{\cos (x^3)},\label{cos x3}
\ee 
satisfy
\be
\frac{D\textrm{Sin }x}{Dx} &=& \textrm{Cos }x,\label{DSin}\\
\frac{D\textrm{Cos }x}{Dx} &=& \ominus\textrm{Sin }x=-\textrm{Sin }x,
\ee
where $\ominus x=0\ominus x =\sqrt[3]{-x^3}=-x$, and
\be
\textrm{Sin}^2x\oplus \textrm{Cos}^2x
&=&
\sqrt[3]{\sin^2(x^3)+\cos^2(x^3)}
=
1.
\ee
$\textrm{Sin }x $ and $\textrm{Cos }x $ are essentially the chirp signals known from signal analysis (Fig.~\ref{Fig1} and Fig.~\ref{Fig2}).

It is instructive to compare (\ref{DSin}) with the Newtonian derivative 
\be
\frac{d \textrm{Sin }x}{dx}
=
\frac{x^2 \cos (x^3)}{\sin ^{\frac{2}{3}}(x^3)},\label{dSin}
\ee
defined with respect to the `standard' arithmetic (Fig.~\ref{Fig2}).

Even more intriguing examples occur if one considers derivatives of functions $A: \mathbb{X}\to \mathbb{Y}$ where the domain and the image of $A$ involve different arithmetics. Let $\mathbb{X}=\mathbb{R}_+$, $\mathbb{Y}=\mathbb{R}$, $f_{\mathbb{X}}(x)=\ln x$, $f_{\mathbb{Y}}(x)=x^3$. The arithmetic operations in $\mathbb{X}$ read explicitly
\be
x_1\oplus_\mathbb{X}x_2
&=&
f_\mathbb{X}^{-1}\big(f_{\mathbb{X}}(x_1)+f_{\mathbb{X}}(x_2)\big)
=
e^{\ln x_1+\ln x_2}=x_1x_2,\\
x_1\ominus_\mathbb{X}x_2
&=&
f_\mathbb{X}^{-1}\big(f_{\mathbb{X}}(x_1)-f_{\mathbb{X}}(x_2)\big)
=
e^{\ln x_1-\ln x_2}=x_1/x_2,\\
x_1\otimes_\mathbb{X}x_2
&=&
f_\mathbb{X}^{-1}\big(f_{\mathbb{X}}(x_1)f_{\mathbb{X}}(x_2)\big)
=
e^{\ln x_1\ln x_2}=x_1^{\ln x_2}=x_2^{\ln x_1},\\
x_1\oslash_\mathbb{X}x_2
&=&
f_\mathbb{X}^{-1}\big(f_{\mathbb{X}}(x_1)/f_{\mathbb{X}}(x_2)\big)
=
e^{\ln x_1/\ln x_2}=x_1^{1/\ln x_2}.
\ee
Neutral elements in $\mathbb{X}$ are given by
\be
1_\mathbb{X}&=&f_\mathbb{X}^{-1}(1)=e^1=e,\\
0_\mathbb{X}&=&f_\mathbb{X}^{-1}(0)=e^0=1.
\ee
A negative of $x\in\mathbb{X}$ is defined as
\be
\ominus_\mathbb{X}x &=& 0_\mathbb{X}\ominus_\mathbb{X}x=f_\mathbb{X}^{-1}\big(-f_\mathbb{X}(x)\big)=e^{-\ln x}=1/x\in\mathbb{R}_+.
\ee
As we can see, numbers negative with respect to the arithmetic from $\mathbb{X}$ are positive if treated in the usual sense.
The unique solution $A:\mathbb{X}\to \mathbb{Y}$ of
\be
\frac{D A(x)}{Dx}=A(x),\quad A(0_\mathbb{X})=1_\mathbb{Y},\label{2exp eq}
\ee
turns out to be
\be
A(x)=f^{-1}_\mathbb{Y}\left(e^{f_\mathbb{X}(x)}\right)=\sqrt[3]{e^{\ln x}}=\sqrt[3]{x}.\label{2exp sol}
\ee
Indeed, first of all, 
\be
A(0_\mathbb{X})
&=&
\sqrt[3]{1}
=1=1_\mathbb{Y}.
\ee
Recalling that multiplication in $\mathbb{Y}$ is unchanged, we check directly from definition (cf. \cite{ACK3}):
\be
\frac{DA(x)}{Dx}
&=& \lim_{h\to 0}\Big(A\big(x\oplus_{\mathbb{X}}f^{-1}_{\mathbb{X}}(h)\big)\ominus_{\mathbb{Y}}A(x)\Big)\oslash_{\mathbb{Y}}f^{-1}_{\mathbb{Y}}(h)
\\
&=& \lim_{h\to 0}\Big(\sqrt[3]{x\oplus_{\mathbb{X}}e^h}\ominus_{\mathbb{Y}}\sqrt[3]{x}\Big)/\sqrt[3]{h}
\\
&=& \lim_{h\to 0}\sqrt[3]{x\frac{e^h-1}{h}}=\sqrt[3]{x}=A(x).
\ee
The exponent satisfies
\be
A(x_1\oplus_{\mathbb{X}}x_2)
=
A(x_1x_2)=\sqrt[3]{x_1x_2}=\sqrt[3]{x_1}\sqrt[3]{x_2}=A(x_1)A(x_2)=A(x_1)\odot_\mathbb{Y}A(x_2),
\ee
as expected. The results are counterintuitive but consistent. The bijection $f_{\mathbb{X}}(x)=\ln x$ is a simplest example of an information channel associated with human or animal nervous system (the Weber-Fechner law; this is why decibels correspond to a logarithmic scale \cite{MC3}).

As final two examples consider first $f_{\mathbb{X}}(x)=x$, $f_{\mathbb{Y}}(x)=\ln x$. The non-Newtonian derivative reads explicitly
\be
\frac{DA(x)}{Dx}
&=&
\lim_{h\to 0}\Big(A(x+h)\ominus_\mathbb{Y} A(x)\Big)\oslash_\mathbb{Y} h_\mathbb{Y}\\
&=&
\lim_{h\to 0}e^{\big(\ln A(x+h)-\ln A(x)\big)/ h}=e^{A'(x)/A(x)}.
\ee
Here $A'(x)=dA(x)/dx$ is the Newtonian derivative.  Let us now solve
\be
\frac{DA(x)}{Dx}
&=&
A(x),\quad A(0)=1_\mathbb{Y}=f^{-1}_\mathbb{Y}(1)=e, 
\ee
an equation equivalent to
\be
e^{A'(x)/A(x)}=A(x).
\ee
By the general formula we know that this must be the non-Newtonian exponent,
\be
A(x)&=&f^{-1}_\mathbb{Y}\left(e^{f_\mathbb{X}(x)}\right)=e^{e^x}.
\ee
Secondly, let $f_{\mathbb{X}}(x)=\ln x=f_{\mathbb{Y}}(x)$. Then 
\be
\frac{DA(x)}{Dx}
&=&
e^{x A'(x)/A(x)}.
\ee
Here values of non-Newtonian and Newtonian exponents coincide,
\be
A(x)&=&f^{-1}_\mathbb{Y}\left(e^{f_\mathbb{X}(x)}\right)=e^{e^{\ln x}}=e^x,
\ee
but their domains are different. Both types of differentiation have been extensively studied in the literature, with numerous applications \cite{mul1,mul2,mul3,mul4,mul5,mul6}. The variety of applications, from signal processing to economics, is not that surprising if one realizes that $\ln x$ represents a neuronal information channel \cite{MC3}. The two non-Newtonian derivatives represent here a perception of change, and not the change itself.

Armed with these intuitions, we are ready to apply the formalism to waves on Koch-type fractals.

\section{Koch curve supported on unit interval}

\begin{figure}
\includegraphics[width=6 cm]{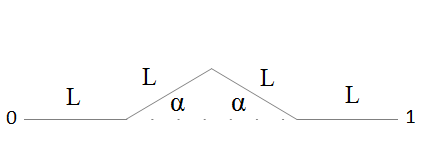}

\includegraphics[width=11cm]{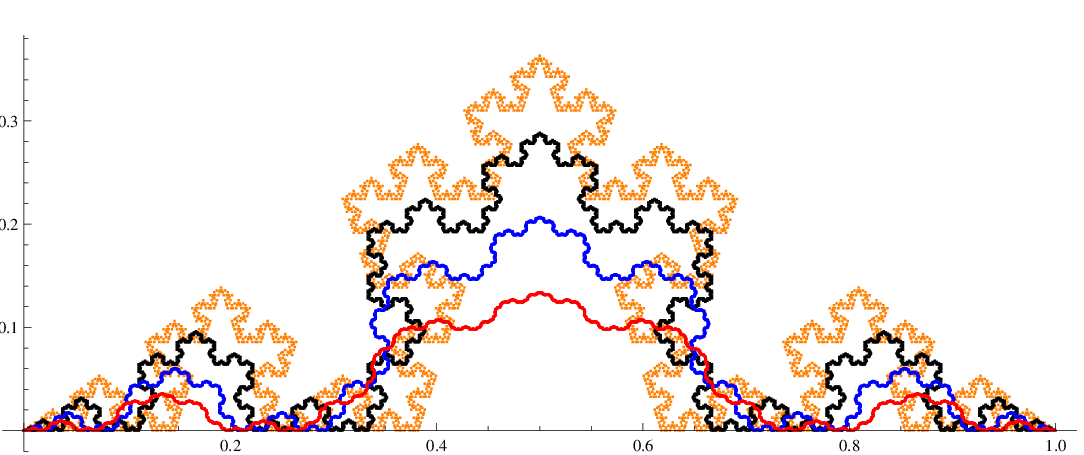}
\caption{Koch curves and their generator (the upper inset) parametrized by $\alpha$ and corresponding to (\ref{0'})-(\ref{z}). From highest to lowest: $\alpha=\pi/2.5$, $\alpha=\pi/3$, $\alpha=\pi/4$, $\alpha=\pi/6$. }
\label{Fig3}
\end{figure}

For convenience we represent $\mathbb{R}^2$ by $\mathbb{C}$. Let us begin with the Koch curve $K_{[0,1]}\subset \mathbb{C}$, beginning at 0 and ending at 1 (Fig.~\ref{Fig3}).
A point $z\in K_{[0,1]}$ can be parametrized by a real number in quaternary representation,
\be
y=(0.q_{1}\dots q_{j}\dots)_4\in [0,1]
\ee
where $q_k=0,1,2,3$. The parametrization is defined by a bijection $g: [0,1]\to K_{[0,1]}$, $z=g(y)$, constructed as follows.
Consider $a=e^{i\alpha}$, $0\leq \alpha\leq \pi/2$, $L=1/(2+2\cos\alpha)$, and
\be
\hat 0(z) &=& Lz,\label{0'}\\
\hat 1(z) &=& L(1+az),\\
\hat 2(z) &=& L(1+a+\bar az),\\
\hat 3(z) &=& L(1+2\cos\alpha +z).\label{3'}
\ee
An $n$-digit point $z\in K_{[0,1]}$ corresponding to $y=(0.q_{1}\dots q_{n})_4$, $q_{n}\neq 0$, is given by
\be
\hat q_{1}\circ \dots \circ\hat q_{n}(0) &=& g(y)\label{z}
\ee
(value at 0 of the composition of maps).
If $y_n=(0.q_{1}\dots q_{n})_4$ is a Cauchy sequence convergent to $y=\lim_{n\to\infty}y_n$, then $g(y)=\lim_{n\to\infty}g(y_n)$.
Curves from Fig.~\ref{Fig3} are the images $g\big([0,1]\big)$ for various $\alpha$. $g$ is one-one, so it defines the inverse bijection $g^{-1}=f: K_{[0,1]}\to [0,1]$.

\begin{figure}
\includegraphics[width=11cm]{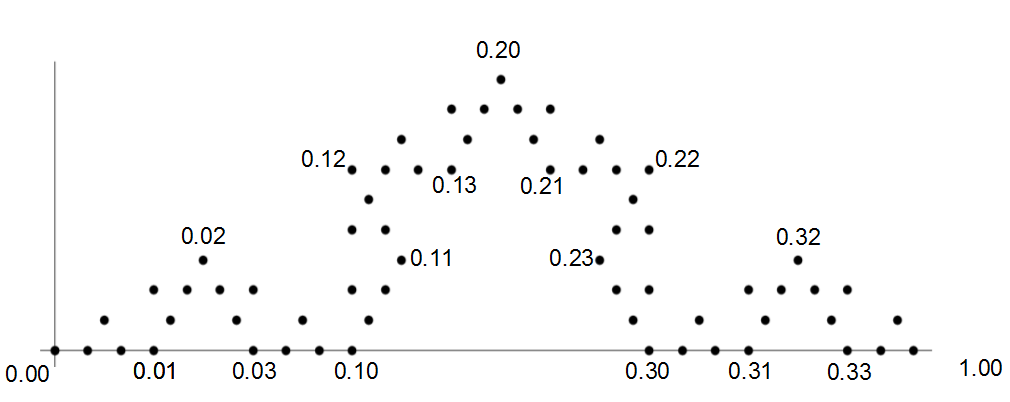}
\caption{Link between a vertex position in the $\alpha=\pi/3$ Koch curve and its numbering by $y=(0.q_{1}\dots q_{j}\dots)_4\in [0,1)$ in quaternary representation. Rescaling the unit segment three times we obtain a smaller copy of the Koch curve. The corresponding vertices of the two curves are numbered by identical digits, with digital separators shifted by one position. The number $y$ can be thus regarded as a Hausdorff measure of the part of the Koch curve extending between the origin and the vertex, if we normalize the measure to 1 on the segment $[0,1)$. The rule applies to all the Koch curves generated by (\ref{0'})--(\ref{3'}).}
\label{numeracja}
\end{figure}

In order to have a better feel of our bijection let us have a look at the relation between the standard $\pi/3$ Koch curve and its quaternary parametrization, as illustrated in Fig.~\ref{numeracja}. Decreasing the initiator \cite{ES,ES1} of the Koch curve three times is equivalent to dividing each vertex number by four (i.e. shifting left the decimal separator by one position). The bijection is therefore equivalent to a parametrization of the Koch curve by its Hausdorff integral, in exact analogy to the construction of Epstein and {\'S}niatycki  \cite{ES,ES1}. The authors of \cite{ES,ES1} begin with the integral and obtain derivatives by means of the fundamental theorem of calculus. The arithmetic approach begins with the derivative, and then the integral is defined through the fundamental theorem of non-Newtonian calculus. 
It is also an appropriate place to mention the notion of the Hausdorff derivative and its generalizations (for a recent review see \cite{He}). In this approach one, roughly speaking, replaces $dx$ by $(dx)^d$ where $d$ is the similarity dimension of a `fractal'. In the Koch example this practically means that having two points separated by a finite distance $\Delta x$ in $\mathbb{R}$ one increases them by performing a single step of the Koch-curve generating algorithm, $\Delta x\mapsto (\Delta x)^d$, replacing initiators by generators (zeroth iteration by first iteration). It makes sense in practical applications to approximate prefractals, such a those occurring in modeling of nanofibers, since it replaces the distance in space by the distance along a fiber, simultaneously assuming that there exists a minimal length beyond which the modeling breaks down. However, if one treats a true Koch fractal as a fiber, the problem is that the length of any of its segments would be infinite, and the limit would not exist. The approaches from \cite{ES,ES1} and the one advocated in the present paper do not suffer from this drawback. Let us also note here that the more traditional approaches to fractal analysis \cite{Kigami,Strichartz} have not managed to formulate any calculus on fractals of a Koch-curve type.

\begin{figure}
\includegraphics{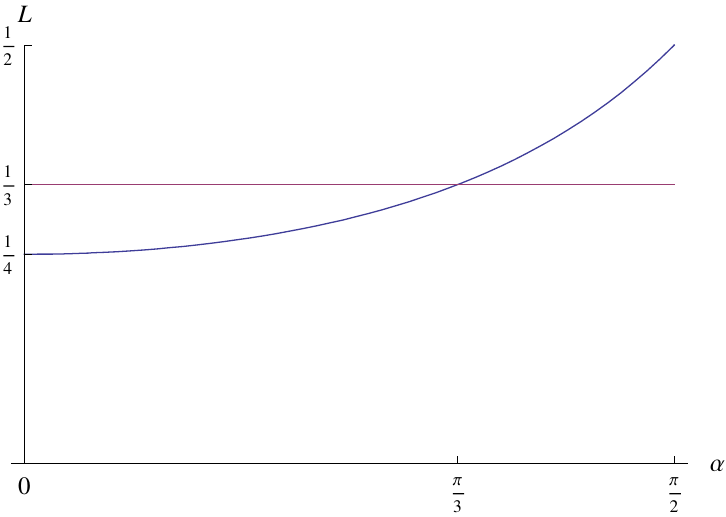}

\includegraphics{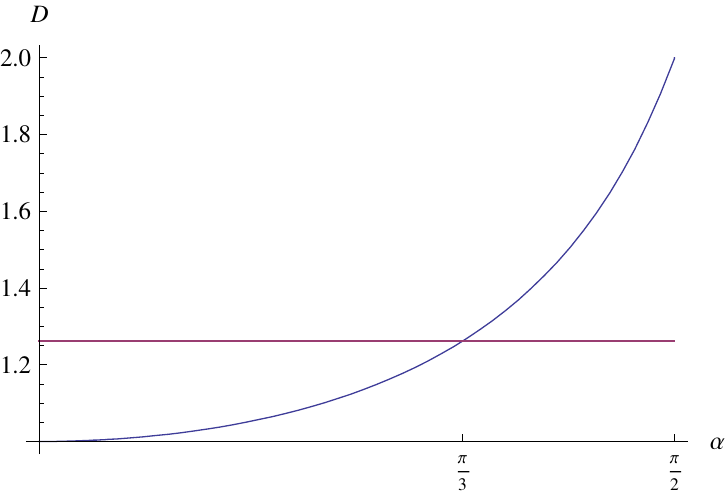}
\caption{Similarity dimension $D$ and the length $L$ of the generator from Fig.~1 as functions of $\alpha$. The horizontal lines show the values for the standard $\pi/3$ Koch curve.}
\label{Fig4}
\end{figure}

For $\alpha=\pi/3$ we obtain the standard curve, generated by equilateral triangles.
Similarity dimension of a curve generated by (\ref{0'})-(\ref{z}) is given by (Fig.~\ref{Fig4})
\be
D &=&\frac{\log 4}{\log (2+2\cos\alpha)}.\label{D}
\ee
There are many ways of extending the Koch curve from $K_{[0,1]}$ to $K_{\mathbb{R}}$. For example, let $K_{[k,k+1]}$ be the curve $K_{[0,1]}$ shifted according to $z\mapsto z+k$, $k\in \mathbb{Z}$. Then $K_{\mathbb{R}}=\cup_{k\in \mathbb{Z}}K_{[k,k+1]}$ is a periodic Koch curve, with the bijection
$f: K_{\mathbb{R}}\to \mathbb{R}$ constructed from appropriately shifted maps $g$ defined above. Non-periodic but self-similar extensions can be obtained by shifts and rescalings. From our point of view the only condition we impose on $f$ is the continuity of $g=f^{-1}$ at 0, i.e. $\lim_{y\to 0_-}g(y)=\lim_{y\to 0_+}g(y)=g(0)$. We take $g(0)=0$.

Combining the generalized Koch curves we can construct a curve which is in a one-one relation with $\mathbb{R}$, with explicitly given bijection $f$, and whose fractal dimensions vary from segment to segment in a prescribed way. This type of generalization may be useful for applications involving realistic coastlines, whose fractal dimensions coincide with the data described by the Richardson law \cite{Mandelbrot}.
In what follows, we will concentrate on the simple case $\alpha=\pi/3$, $L=1/3$, of the standard Koch curve.

\section{Wave equation on Koch curves}

First of all, let us assume we discuss a real-valued field, whose evolution on the Koch curve $\mathbb{X}=K_{\mathbb{R}}$ is described with respect to a `normal' non-fractal time $t$. The field is thus represented by $\mathbb{R}\times \mathbb{X}\mapsto \Phi_t(x)\in \mathbb{R}$, with $x\in \mathbb{X}$. Since $\mathbb{Y}=\mathbb{R}$ we take $f_\mathbb{Y}=\textrm{id}_\mathbb{R}$. (Although $f_\mathbb{Y}(y)=y^3$ or any other bijection would do as well, leading to a different behavior of the wave.) The wave equation is
\be
\frac{1}{c^2}\frac{d^2}{dt^2}\Phi_t(x)-\frac{D^2}{Dx^2}\Phi_t(x)=0,
\ee
where
\be
\frac{d}{dt}\Phi_t(x) &=& \lim_{h\to 0}\Big(\Phi_{t+h}(x)-\Phi_t(x)\Big)/h,\\
\frac{D}{Dx}\Phi_t(x) &=& \lim_{h\to 0}\Big(\Phi_t\big(x\oplus_\mathbb{X} f^{-1}_\mathbb{X}(h)\big)-\Phi_t(x)\Big)/h.
\ee
We search solutions in the form (here $y=ct$)
\be
\Phi_t(x)
=
A(x,y)+B(x,y),
\ee
where
\be
\left(\frac{d}{dy}-\frac{D}{Dx}\right)A(x,y)=
\left(\frac{d}{dy}+\frac{D}{Dx}\right)B(x,y)\equiv 0,
\ee
suggesting simply
\be
A(x,y) &=& a\big(f_\mathbb{X}(x)+y\big),\\
B(x,y) &=& b\big(f_\mathbb{X}(x)-y\big),\label{b}
\ee
for some twice differentiable $a,b:\mathbb{R}\to \mathbb{R}$.

\begin{figure}
\includegraphics{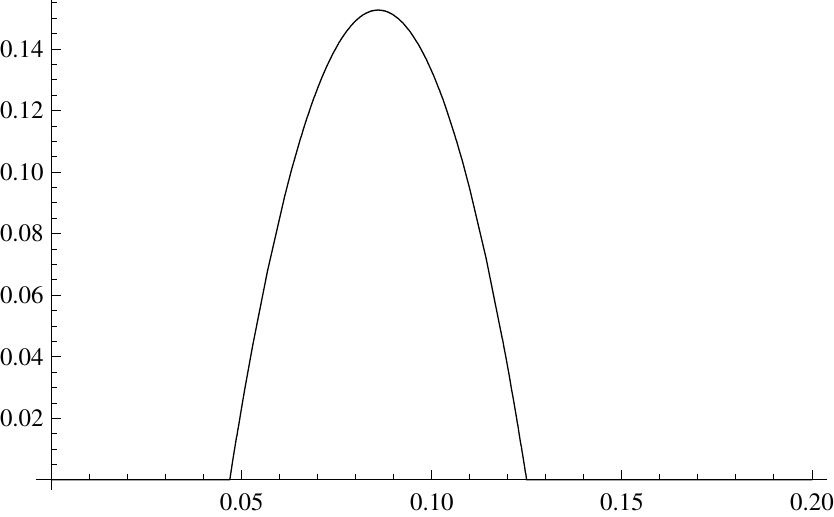}

\includegraphics[width=14cm]{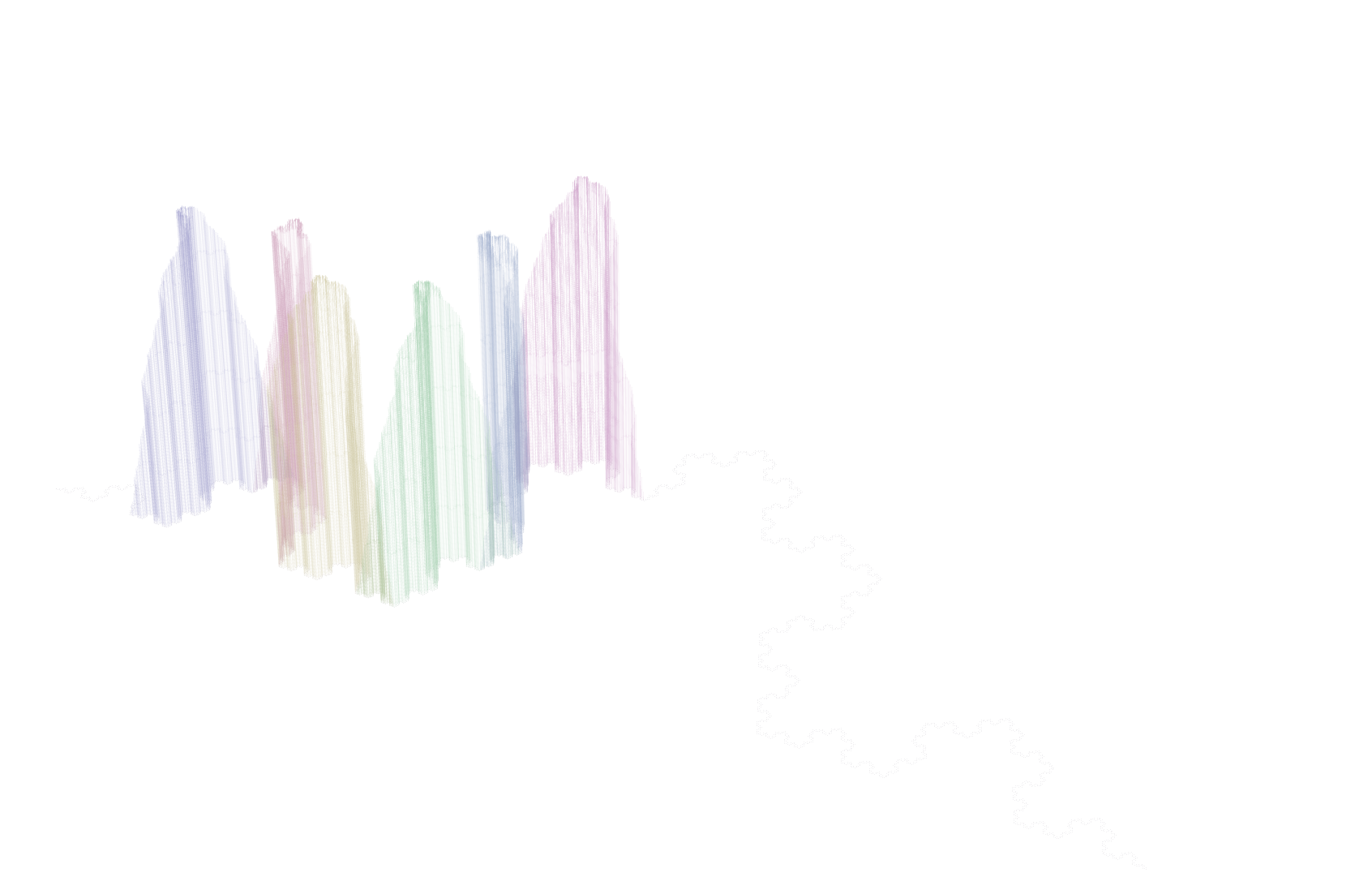}
\caption{`Aurora borealis wave': Six snapshots of $\Phi_t(x)$ propagating to the right along the Koch curve. The upper plot shows the corresponding function $b$ occurring in (\ref{b}).}
\label{Fig5}
\end{figure}

Indeed, from definitions
\be
\frac{D}{Dx}A(x,y)
&=&
\lim_{h\to 0}\frac{A\big(x\oplus_\mathbb{X} f^{-1}_\mathbb{X}(h),y\big)-A(x,y)}{h}
\nonumber\\
&=&
\lim_{h\to 0}\frac{a\Big(f_\mathbb{X}(x)+h+y\big)-a\big(f_\mathbb{X}(x)+y\big)}{h}
\nonumber\\
&\equiv&
\frac{d}{dy}a\big(f_\mathbb{X}(x)+y\big)=\frac{d}{dy}A(x,y).
\ee
One similarly verifies that $d/dy$ and $D/Dx$ commute, and
\be
\frac{D}{Dx}B(x,y)
&\equiv&
-\frac{d}{dy}B(x,y).
\ee
Fig.~\ref{Fig5} shows the dynamics of $\Phi_t(x)$ with $a=0$.
The energy of the wave is given by
\be
E
&=&
\frac{1}{2}\int_{f_\mathbb{X}^{-1}(-\infty)}^{f_\mathbb{X}^{-1}(\infty)}
\left(
\frac{1}{c^2}\left|\frac{d\Phi_t(x)}{dt}\right|^2
+
\left|\frac{D\Phi_t(x)}{Dx}\right|^2
\right)
Dx,
\nonumber\\
\ee
where the integral is defined by (\ref{int}).

Let us explicitly check the time independence of $E$ for the particular case of
$\Phi_t(x)=a\big(f_\mathbb{X}(x)+ct\big)$. Let $a'(x)=da(x)/dx$ be the Newtonian derivative. Then,
\be
E
&=&
\int_{-\infty}^{\infty}
|a'\big(f_\mathbb{X}\circ f_\mathbb{X}^{-1}(x)+ct\big)|^2
dx\\
&=&
\int_{-\infty}^{\infty}
|a'(x)|^2
dx
\ee
is independent of time, as it should be.

\section{Lorentz covariance}

In our model space-time consists of points
$(x^0,x^1) = (ct,x)\in\mathbb{R}\times \mathbb{X}$,
with $\big(x^0,f_{\mathbb{X}}(x^1)\big)\in\mathbb{R}^2$. A Lorentz transformation $x'={\cal L}(x)$, ${\cal L}:\mathbb{R}\times \mathbb{X}\to \mathbb{R}\times \mathbb{X}$, is defined by
\be
\left(
\begin{array}{c}
x'^0\\
x'^1
\end{array}
\right)
&=&
\left(
\begin{array}{c}
L{^0}{_0} x^0+L{^0}{_1}f_{\mathbb{X}}(x^1)\\
f_{\mathbb{X}}^{-1}\Big(L{^1}{_0} x^0+L{^1}{_1}f_{\mathbb{X}}(x^1)\Big)
\end{array}
\right),\label{Lor}
\ee
or, equivalently, by
\be
\left(
\begin{array}{c}
x'^0\\
f_{\mathbb{X}}(x'^1)
\end{array}
\right)
&=&
\left(
\begin{array}{cc}
L{^0}{_0} & L{^0}{_1}\\
L{^1}{_0} & L{^1}{_1}
\end{array}
\right)
\left(
\begin{array}{c}
x^0\\
f_{\mathbb{X}}(x^1)
\end{array}
\right)\label{Lor0},
\ee
where $L\in {\rm SO}(1,1)$. (\ref{Lor}) implements a nonlinear action of the group ${\rm SO}(1,1)$, and reduces to the usual representation if $\mathbb{X}=\mathbb{R}$ and $f_{\mathbb{X}}(x^1)=x^1$. Transformations (\ref{Lor}) form a group.

In order to prove Lorentz invariance of the wave equation let us first note that its solution
\be
\Phi_t(x) &=&
a\big(f_\mathbb{X}(x^1)+x^0\big)+b\big(f_\mathbb{X}(x^1)-x^0\big)
\nonumber\\
&=&
\phi\big(x^0,f_\mathbb{X}(x^1)\big),
\ee
defines a function $\phi$, satisfying (due to triviality of $f_{\mathbb{Y}}$)
\be
\frac{D\Phi_t(x)}{Dx}
&=&
\frac{\partial\phi\big(x^0,f_\mathbb{X}(x^1)\big)}{\partial f_\mathbb{X}(x^1)},\\
\frac{1}{c}\frac{d\Phi_t(x)}{dt}
&=&
\frac{\partial\phi\big(x^0,f_\mathbb{X}(x^1)\big)}{\partial x^0}.
\ee
Accordingly, the wave equation takes the standard form
\be
\left(
\frac{1}{c^2}\frac{\partial^2}{\partial t^2}
-
\frac{\partial^2}{\partial f_\mathbb{X}(x^1){}^2}
\right)\phi\big(x^0,f_\mathbb{X}(x^1)\big)
=0.\label{standard}
\ee
It is invariant under (\ref{Lor0}) if $\phi$ transforms by
\be
\phi'\big(x'^0,f_\mathbb{X}(x'^1)\big)
&=&
\phi\big(x^0,f_\mathbb{X}(x^1)\big),
\ee
which is equivalent to the scalar-field transformation $\Phi'_{t'}(x')=\Phi_t(x)$.

Replacing $\mathbb{R}\times \mathbb{X}$ by a more general case $\mathbb{X}_0\times \mathbb{X}_1$, $f_{\mathbb{X}_j}:\mathbb{X}_j\to \mathbb{R}$, one arrives at a Lorentz invariant wave equation (with both space-time derivatives appropriately defined), and Lorentz transformations
\be
\left(
\begin{array}{c}
x'^0\\
x'^1
\end{array}
\right)
&=&
\left(
\begin{array}{c}
f_{\mathbb{X}_0}^{-1}\Big(L{^0}{_0} f_{\mathbb{X}_0}(x^0)+L{^0}{_1}f_{\mathbb{X}_1}(x^1)\Big)\\
f_{\mathbb{X}_1}^{-1}\Big(L{^1}{_0} f_{\mathbb{X}_0}(x^0)+L{^1}{_1}f_{\mathbb{X}_1}(x^1)\Big)
\end{array}
\right).\label{Lor-gen}
\ee
A generalization to space-times constructed by Cartesian products of arbitrary numbers of fractals is now obvious.

\section{Conclusions}

To conclude, we have obtained a wave that propagates along a Koch-type curve. The wave possesses finite conserved energy and satisfies the usual wave equation, formulated with respect to appropriately defined derivatives. The derivatives are not the ones we know from the standard mathematical literature of the subject, but are very natural and easy to work with. The solution we have found is the general one, a fact following from the standard form (\ref{standard}) of the wave equation. The velocity of the wave is intriguing. On the one hand, it is described by the parameter $c$ in the wave equation. On the other hand, however, the length of any piece of a fractal coast is infinite and yet the wave moves from point to point in a finite time, and with speed that {\it looks\/} finite and natural. This is possible since the fractal sum $z=x\oplus_{\mathbb{X}} y$ of two points in a Koch curve is uniquely defined in spite of the apparently `infinite' distances between $x$, $y$, $z$ and the origin 0. Another interesting aspect of the resulting motion is the lack of difficulties with combining non-fractal time with fractal space. Lorentz transformations in the corresponding space-time have been constructed, and Lorentz invariance of the wave equation has been proved. Fractal arithmetic automatically tames the infinities inherent in the length of the curve. It would not be very surprising if our fractal calculus found applications also in other branches of physics, where finite physical results are buried in apparently infinite theoretical predictions.

\end{document}